\documentclass[11pt]{article}
\usepackage{amsfonts}
\usepackage{setspace}

\newcommand{\triangledown}       {\nabla}
\newcommand{\mti}       {m\rightarrow\infty}
\newcommand{\N}       {{\cal N}}
\renewcommand{\theequation}{\thesection.\arabic{equation}}
\renewcommand{\Re}{I\!\!R}
\renewcommand{\amalg}{1\!\!1}
\def\qmq#1{\quad\mbox{#1}\quad}

\newcommand{\ua}       {\mbox{\boldmath$a$}} 
 
\newcommand{\ub}       {\mbox{\boldmath$b$}}

\newcommand{\ue}       {\mbox{\boldmath$e$}}

\newcommand{\uV}       {\mbox{\boldmath$V$}}

\newcommand{\ux}       {\mbox{\boldmath$x$}}

\newcommand{\ubeta}             {\mbox{\boldmath$\beta$}}
\newcommand{\ugamma}            {\mbox{\boldmath$\gamma$}}
\newcommand{\udelta}            {\mbox{\boldmath$\delta$}}

\newcommand{\uiota}             {\mbox{\boldmath$\uiota$}}

\newcommand{\uSigma}            {\mbox{\boldmath$\Sigma$}}

\newcommand{\be}               {\begin{equation}}
\newcommand{\ee}               {\end{equation}}
\newcommand{\bea}              {\begin{eqnarray}}
\newcommand{\eea}              {\end{eqnarray}}
\newcommand{\ba}               {\begin{array}}
\newcommand{\ea}               {\end{array}}
\newcommand{\nn}               {\nonumber}

\begin{document}

\begin{center}
\LARGE
{\bf Nonnegative  mean squared prediction error estimation  
 in small area estimation}
\vskip.2in

\normalsize

{\sc By} SOUMENDRA N. LAHIRI {\sc And} TAPABRATA MAITI\\

{\it Department of Statistics, Iowa State University, Ames, Iowa, USA}

\vspace{.3in}

{\sc Summary}
\end{center}

Small area estimation has received  enormous attention in recent years 
due to its wide range of application, particularly 
in policy making decisions. The variance based on 
direct sample size of small area estimator is unduly 
large and there is a need of constructing model based 
estimator with low mean squared prediction error (MSPE).  Estimation 
of MSPE and in particular the bias correction of MSPE 
plays the central piece of small area estimation research. In this 
article, a new technique of bias correction 
for  the estimated MSPE is proposed. It is shown that that the 
new MSPE estimator attains the same level of bias correction 
as the existing estimators based on straight Taylor expansion 
and jackknife methods. However, unlike the existing methods, the 
proposed estimate of MSPE is always nonnegative. Furthermore,
the proposed method can be used for  general two-level small
area models where the variables at each level can be discrete
or continuous and, in particular,  be nonnormal.   
 \vspace*{2ex}

\noindent
{\it Some key words:} Best predictor; Bootstrap; Mean squared prediction 
error; small area. 

\section{\sc Introduction}

Small area estimation  is an important statistical research area due
 to its growing demand from public and private agencies. The variance 
of a  small area estimator is unduly large due to smallness of the
area-level  sample size. Use of models has proven to be unavoidable 
to control  the mean squared prediction error (MSPE) of 
a small area predictor.  The bias correction of the 
 estimated MSPE is the central piece of small
 area estimation research. See Rao (2003) and references therein 
for  a detailed  discussion. 

The standard small area models are usually two-level models,
where  one is a sampling  model and the other one
 is a population model. Prasad and Rao (1990) 
assumed normality at  both levels and used ANOVA estimates 
of the  model  parameters to derive  second order correct MSPE 
estimates. Lahiri  and Rao (1995) relaxed the normality
assumption at the population level and 
re-establish the  Prasad-Rao result on second order correct MSPE 
estimation. Datta and Lahiri (2000) investigated properties of  
Prasad-Rao (PR)  type MSPE estimators for maximum likelihood and 
restricted maximum likelihood estimates of the model parameters,
retaining the   normal distribution assumption at 
both the levels. Recently, Jiang, Lahiri 
and Wan (2002) proposed a jackknife based MSPE estimators where 
normality is not a requirement. However,  the Jiang-Lahiri-Wan
 method (JLW) requires a closed form expression for the  posterior risk 
which is not often  available (e.g., the  binomial-normal model). 
Moreover, the JLW estimator  has the undesirable property that 
it may produce  negative MSPE  estimates (Bell, 2002). Although 
the PR type MSPE estimates are nonnegative for the 
normal-normal case, the nonnegativity property 
 is unknown for other situations. 
 The PR type MSPE estimators correct the bias 
of the estimated MSPE using Taylor's expansion. 
On the other hand, the  JLW MSPE estimator corrects the bias of
the  posterior risk using the jackknifing method. In this article,
 we propose a new technique of MSPE bias 
correction which attains the same level of accuracy  
as that of the PR type or the  JLW MSPE estimators. In addition,
 the new MSPE estimates are  guaranteed to be nonnegative. Moreover, 
the new  method is valid for any family of  parametric distributions, 
discrete or continuous. Thus,  unlike the traditional 
methods, neither the normality assumption nor the choice of a specific 
parameter estimation  method are required  for the validity
of the  proposed approach.

The organization of the paper is as follows: The next section introduces
the two-level small area models and  discusses the 
existing MSPE estimation methods in  this  framework.
 The Section 3  proposes the new MSPE estimator. 
Some technical properties of the proposed estimator
 are discussed and  compared  with the existing methods 
in Section 4. Section 5 reports  finite sample properties of the new 
estimator  using a simulation study. Some  conclusions and comments are
 made in  Section 6. Proofs of the technical results are given in the 
Appendix.

\section{\sc Existing Methods of MSPE Estimation}
\setcounter{section}{2}
\setcounter{equation}{0}
Consider the two-level small area model
\bea
y_i & = & \theta_i +e_i, \; e_i \sim F(.;D_i) \\
\theta_i & = & \ux_i^T\ubeta +u_i, \; u_i \sim G(.;\ugamma)
\eea
$i=1,\cdots,m$, where $y_1,\cdots,y_m$ are direct estimators 
with sampling errors $e_1,\cdots,e_m$, independently distributed 
with cumulative distribution functions $F(.;D_1),\cdots,F(.;D_m)$,
respectively; $u_1,\cdots,u_m$ are independent 
and identically distributed (iid) 
random variables with common distribution function $G(.;\ugamma)$, 
and $\ux_1,\cdots\ux_m$ are $p$-dimensional nonrandom covariates. 
We suppose that the unknown parameters of the model are given by 
the regression parameter $\ubeta$ and the $p$- dimensional 
parameter $\ugamma$ of the random effects distribution $G(.;\ugamma)$ 
in (2.2), but the values of $D_1,\cdots D_m$ are known, as typically 
assumed in two-level small area models. We also assume that the 
sampling errors $e_1,\cdots,e_m$ and random effects $u_1,\cdots,u_m$ 
are mutually independent with $E(e_i)=0=E(u_i),i=1,\cdots,m$.  
Note that neither the $e_i$'s nor the $u_i$'s 
are required to be normally distributed. In fact, under (2.1) 
and (2.2), the $e_i$'s and the $u_i$'s are allowed to have 
arbitrary parametric families of  discrete or continuous 
distributions.

Suppose that the quantity of interest for prediction is given by 
\be
h(\theta_i), i=1,\cdots,m
\ee
for some smooth function $h: \Re\rightarrow \Re$. For example, 
$h(x)=x$ is the most commonly used function, which may correspond 
to area level means or totals. An important example of $h(.)$ 
includes exponentiation in U.S. Census Bureau's ongoing Small Area Income 
and Poverty Estimation 
(SAIPE) project. For county level poverty estimation in SAIPE, 
the model (2.1) and (2.2) applies after log transformation 
of the original data.

The best predictor (BP) of $h(\theta_i)$ is given by 
\be
H_i(\udelta)\equiv E_{\udelta}(h(\theta_i)|y_i), ~~i=1,\cdots, m
\ee
where $\udelta=(\ubeta^T,\ugamma^T)^T$ is the vector of model 
parameters. Since the true value of $\udelta$ is unknown, 
$H_i(\udelta)$ is not directly usable in practice. 
It is customary to substitute an estimator $\hat{\udelta}$, 
say, of $\udelta$ and predict $h(\theta_i)$ by using the
 estimated best predictor (EBP) as 
\be
H_i(\hat{\udelta}), i=1,\cdots,m.
\ee
Performance of the EBP is measured by the mean squared 
prediction error (MSPE):
\be
M_i(\udelta)=E_{\udelta}\left(H_i(\hat{\udelta})-
h(\theta_i)\right)^2, i=1,\cdots, m.
\ee
Further, like the EBP, an estimator of the MSPE is obtained by 
$M_i(\hat{\udelta}),\; i=1,\cdots,m.$ However, as pointed out 
by Prasad and Rao (1990) in their seminal paper, this naive 
plug-in estimator is not very useful. 
To appreciate why, note that $M_i(\udelta)$ can be  decomposed as 
\bea
M_i(\udelta) & = & E_{\udelta}\left(H_i(\hat{\udelta})
-h(\theta_i)\right)^2 
\nn\\ 
 & = & E_{\udelta}\left(H_i(\udelta)-h(\theta_i)\right)^2+
E_{\udelta}\left(H_i(\hat{\udelta})-H_i(\udelta)\right)^2  \nn\\
 & \equiv & M_{1i}(\udelta)+M_{2i}(\udelta), i=1,\cdots,m,
\eea
where the cross-product term vanishes as a consequence of
 the fact that 
$E_{\udelta}\left(H_i(\udelta)-h(\theta_i)\right)Z=0$ 
for any $\sigma\langle y_1,\cdots,y_m\rangle$-measurable random 
variable $Z$. 
In (2.7), the first term $M_{1i}(\udelta)$ is the optimal 
prediction error using the unknown ideal predictor  $H_i(\udelta)$ 
and is of order $O(1)$ as $m\rightarrow \infty$. The second term 
$M_{2i}(\udelta)$ arises from the error in  estimating the 
unknown model parameters  $\udelta$ in the BP $H_i(\udelta)$, 
and, typically, it is of the order $O(m^{-1})$ as 
$m\rightarrow \infty$.  Prasad and Rao (1990) 
showed that by substituting $\hat{\udelta}$ for $\udelta$ to define 
the naive plug-in estimator 
$$
M_i(\hat{\udelta})=M_{1i}(\hat{\udelta})
+M_{2i}(\hat{\udelta}), 
$$
one introduces an additional bias  of 
the order $O(m^{-1})$, which is of the same order as the second 
term $M_{2i}(\udelta)$ in (2.7). As a  result, the naive 
estimator has  a masking effect on the bias of the 
EBP and hence, it  is not a good estimator of $M_i(\udelta)$,
particularly when $m$ is not too large.

Prasad and Rao (1990) suggested a bias corrected estimator of 
the MSPE  $M_i(\udelta)$ for a normal-normal model. The key idea 
there is to estimate the (leading term of the ) bias of 
$M_i(\hat{\udelta})$ using explicit analytical  expressions. 
The bias corrected estimated MSPE, proposed by Prasad  and Rao (1990), 
 is of the form 
\be
\hat{M}_i^{PR}=M_{1i}(\hat{\udelta})-
\widehat{Bias}_i^{PR}+M_{2i}(\hat{\udelta}),
\ee
where $\widehat{Bias}_i^{PR}$ is obtained by estimating higher 
order terms in the Taylor's expansion of the function $M_{1i}(.)$ 
around $\udelta$. 

An alternative approach, put forward by Jiang, Lahiri and Wan (2002), 
involves using the jackknife method to correct the $O(m^{-1})$-order 
bias term in the naive estimator $M_i(\hat{\udelta})$. More 
specifically,  the bias-corrected estimator of $M_i(\udelta)$ of 
JLW is given by 
\be
\hat{M}_i^{JLW}=M_{1i}(\hat{\udelta})-\widehat{Bias}_i^{JLW}+M_{2i}
(\hat{\udelta}),
\ee
where $\widehat{Bias}_i^{JLW}$ is the Jackknife estimator of the bias of 
$M_{1i}(\hat{\udelta})$. 

Although, the estimators $\hat{M}_i^{PR}$ and  $\hat{M}_i^{JLW}$ have 
superior bias properties, an undesirable feature of both of these 
estimators is that they may produce negative MSPE estimate with 
positive probabilities. This results from the sampling variability 
of the bias estimators, which may dominate the value of the 
unadjusted naive estimator $M_i(\hat{\udelta})$
and thereby, may lead to a negative value of the bias corrected 
MSPE estimators. 

In this article, we propose a different approach to bias 
correction that is guaranteed to produce a nonnegative 
estimate of the MSPE. The key idea here is to tilt suitably  the 
value of $\hat{\udelta}$, an  initial estimator of $\udelta$, 
{\it before} evaluating the function $M_i(.)$, such that the 
difference between the true MSPE $M_i(\udelta)$ and the value
 of the function  $M_i(.)$ at the new value of  the argument, 
say $\tilde{\udelta}$, is smaller on the average. Since the MSPE 
function $M_i(.)$ is always nonnegative,  the resulting estimator
 of the true MSPE is always nonnegative. The tilted 
value $\tilde{\udelta}$ is constructed from $\hat{\delta}$ using the 
data-values only and hence it is itself an estimator of $\udelta$. 
In constructing $\tilde{\udelta}$, we implicitly correct the bias of 
$M_{1i}(\hat{\udelta})$, by making use of estimates of linear combination 
of the bias and the variance of the initial estimator $\hat{\udelta}$. Here, 
we employ the bootstrap method (Efron, 1979) to derive the bias and variance 
of the estimators of model parameters, although other methods such as 
the jackknife and the delta methods, are equally applicable. 
The details of the correct construction  are given in the next section.  

\section{\sc The Proposed Estimator of the MSPE}
\setcounter{section}{3}
\setcounter{equation}{0}
\subsection{\it Motivation}
To motivate the definition of the  proposed MSPE estimator, consider
 a related deterministic approximation problem, where we wish to 
approximate the value of a smooth function $f: \Re\rightarrow 
\Re$ at a point $a \in \Re$ using its values over an 
interval $I$ containing $a$. For a given $c\neq 0$, 
setting  $x_m\equiv a+\frac{c}{\sqrt{m}}$, $m\geq 1$ and
using  Taylor's expansion,  we get
\be
f(x_m)=f(a)+(x_m-a)f'(a)+\frac 12 (x_m-a)^2f''(a)+O(m^{-\frac 32}).
\ee
This suggests that starting with $x_m$, we may now construct 
a new point $\tilde{x}_m \in I_m$ of the form $\tilde{x}_m=x_m+c_m$, 
such that 
\be
f(\tilde{x}_m)=f(a)+ O(m^{-\frac 32}).
\ee
Indeed, by Taylor's expansion of $f(\tilde{x}_m)$ around $a$,
 we have 
\be
f(\tilde{x}_m) = f(a)+(x_m+c_m-a)f'(a)
+\frac 12 (x_m+c_m-a)^2f''(a)+O(m^{-\frac 32}),
\ee
which satisfies (3.2) if 
\be
(x_m+c_m-a)f'(a)+\frac 12(x_m+c_m-a)^2f''(a)=0.
\ee
Now equation (3.4) can be solved for $c_m$ 
(yielding the solution
$c_m =-\frac{2f'(a)}{f''(a)} - (x_m - a)$)
to find the desired point $\tilde{x}_m$. In deriving the proposed MSPE 
estimator, we employ an extension of this 
  simple idea to the function $f(\cdot)=
M_{1i}(\cdot)$ which is now a function (of several real variables) 
from $\Re^k \rightarrow \Re$. The role of the point $x_m$ 
is played by an initial estimator 
$\hat{\udelta}$. Some additional care is needed to ensure 
that the analog of the tilted point $\tilde{x}_m$, now denoted 
by $\tilde{\udelta}$, is truly an estimator, i.e., a function
 of the data alone and does not involve any parameters 
(e.g., it may not involve  the point ``$a$" in $\tilde{x}_m$,
 which represents the true parameter value $\udelta$ in 
our application).  

\subsection{\it Definition of the proposed estimator}
Let $\hat{\udelta}$ be a given estimator of $\udelta$ and let 
$\ub=\ub(\udelta)=E_{\udelta}(\hat{\udelta}-\udelta)$ denote 
the bias and  $V=V(\udelta)=Var_{\udelta}(\hat{\udelta})$ 
denote the variance matrix of $\hat{\udelta}$ at $\udelta$. 
We shall suppose that   
 some consistent estimators
$\hat{\ub}$ and 
 $\hat{V}$ of the bias and the 
variance matrix of the initial estimator $\hat{\udelta}$
are available.
For example, these may be generated by a suitable
 resampling  method; see Section 5 where we use a 
parametric bootstrap method for this purpose.
To define the tilted estimator of $\udelta$, 
we also suppose that for $i=1,\cdots,m$ 
\be
 \sum_{j=1}^k|M_{1i}^{(j)}(\udelta)| \neq 0, 
\ee
where 
 for a  differentiable function $f: \Re^k\rightarrow \Re$, 
$f^{(j)}$ and $f^{(jl)}$  denote the first and the second 
order partial derivatives with respect to the $j$-th 
co-ordinate and the $(j,l)$-th co-ordinates, respectively, 
$j,l=1,\cdots,k.$ 
Condition (3.5) says that  at least one of the first order 
partial derivatives of the function $M_{1i}(\cdot)$ is nonzero 
at the true value of the parameter $\udelta$ for each $i$. 
For notational simplicity, without loss of 
generality, we suppose that $M_{1i}^{(1)}(\udelta)\neq 0$. Then,
 we define the {\it preliminary-tilted-estimator} of $\udelta$ for 
the $i$-th small area by 
\bea
\bar{\udelta}_i & = & \hat{\udelta}-\left[\sum_{j=1}^kM_{1i}^{(j)}
(\hat{\udelta})\hat{\ub}(j)+\frac 12 \sum_{j=1}^k\sum_{l=1}^kM_{1i}^{(jl)}
(\hat{\udelta})\hat{V}(j,l)\right]\left\{M_{1i}^{(1)}(\hat{\udelta})
\right\}^{-1}\ue_1
\eea
where $\ue_1=(1,0,\cdots,0)^T \in \Re^k$, $\hat{\ub}(j)$ denote 
the $j$-th component of $\hat{\ub}$ and $\hat{V}(j,l)$ denote 
the $(j,l)$-th element of $\hat{V}$. Thus, the estimator
$
\bar{\udelta}_i
$
is obtained from the initial estimator $\hat{\udelta}$
by adding a correction factor to the first component of
$\hat{\udelta}$ only. Note that if instead of 
$M_{1i}^{(1)}(\udelta)$, a different partial derivative 
$M_{1i}^{(l)}(\udelta)$ were nonzero, then we would
define the  preliminary tilted  estimator 
$\bar{\udelta}_i$ by replacing the factor 
$\left\{M_{1i}^{(1)}(\hat{\udelta})
\right\}^{-1}\ue_1$ in (3.6) with 
$\left\{M_{1i}^{(l)}(\hat{\udelta})
\right\}^{-1}\ue_l$, 
where the vector $\ue_l \in \Re^k$
has $1$ in the $l$-th position and zeros elsewhere,
$1\leq l\leq k$.

Next, let $\Delta$ denote the set of possible values of the parameter 
$\udelta$ under the model (2.1) and (2.2). Then the tilted 
estimator of $\udelta$ for the $i$-th small area is defined by 
\be
\tilde{\delta}_i = \left\{\ba{ll}
        \bar{\udelta}_i & \mbox{ if } \bar{\udelta}_i \in \Delta 
\mbox{ and } |M_{1i}^{(1)}(\hat{\udelta})|^{-1} \le (1+\log m)^2 \\
    \hat{\udelta} & \mbox{ otherwise } \ea \right.
\ee
$i=1,\cdots, m$. Thus, if the preliminary estimator $\bar{\udelta}_i$ 
takes values inside the parameter space $\Delta$ and the 
value of the partial derivative $M_{1i}^{(1)}(\hat{\udelta})$
at $\hat{\udelta}$ is not too small, the tilted 
estimator of $\udelta$ is given by $\bar{\udelta}_i$ itself. 
However, in the event that either $\bar{\udelta}_i$ falls outside 
$\Delta$ or $M_{1i}^{(1)}(\hat{\udelta})$ becomes too small, 
we replace it with the original estimator 
$\hat{\udelta}$. Small values of $M_{1i}^{(1)}(\hat{\udelta})$
make the estimator $\bar{\udelta}_i$ unstable and hence,
truncated below. It will be shown in Section 4 that under 
appropriate regularity conditions, the probability of 
getting a preliminary estimator $\bar{\udelta}_i$ 
outside $\Delta$ or that of getting a value of
$M_{1i}^{(1)}(\hat{\udelta})$ below the threshold
$(1+\log m)^{-2}$ tends to zero rapidly 
 as $m\rightarrow \infty$, uniformly in $i$. As a 
consequence, the tilted 
estimator $\tilde{\udelta}_i$ coincides with the 
preliminary tilted estimator $\bar{\udelta}_i$ 
with high probability. The proposed estimator of 
the MSPE is now defined as 
\be
\widehat{M_{i}(\udelta)} = M_{1i}(\tilde{\udelta})+M_{2i}(\hat{\udelta}), 
i=1,\cdots, m. 
\ee
Note that by the construction, the MSPE estimator is 
always positive. In the next section, we show that under 
some regularity conditions, it has a bias that is of the 
order $o(m^{-1})$.
Therefore, the proposed estimator attains 
the same level of accuracy as the previously proposed 
estimators $\hat{M}_i^{PR}$ and $\hat{M}_i^{JLW}$, while at 
the same time, guarantees positivity.

\section{\sc  Theoretical Properties of the Proposed Estimators}
\setcounter{section}{4}
\setcounter{equation}{0}

In this section, we describe some theoretical properties of 
the tilted estimator $\tilde{\udelta}_i$ of (3.7)
and  of the bias 
corrected MSPE estimator $\widehat{M_i(\udelta)}$
of (3.8). For proving the result of this section, we shall
assume the  following regularity conditions
on the model (2.1) and (2.2).\\[.1in]
{\bf Condition S:} 
\begin{enumerate}
\item $\udelta$, the true value of the parameter, is an 
interior point of $\Delta$. 
\item  $M_{1i}$ is twice continuously differentiable on 
$\Delta$ and there exists a constant $C_1\in (0,\infty)$ 
such that 
$$|M_{1i}^{(j)}(\ux)|+|M_{1i}^{(jl)}(\ux)| <C_1 $$
for all $\ux \in \Delta, j,l=1,\cdots,k$ and $i=1,\cdots,m, m\geq 1$.
\item 
(i) $M_{2i}$ is differentiable on $\Delta$.\\
(ii) There exist $C_2, \epsilon_0  \in (0, \infty)$ and $\gamma \in (0,1]$ 
such that 
$$
|M_{1i}^{(jl)}(\ux)-M_{1i}^{(jl)}(\udelta)|+m|M_{2i}^{(j)}
(\ux)-M_{2i}^{(j)}(\udelta)|< C_2\|\ux-\udelta\|^{\gamma}
$$
for all $\ux \in \N\equiv \{\|\ux-\udelta\|\leq \epsilon_0\}$ 
for $j,l=1,\cdots,k; i=1,\cdots,m, m\geq 1$. \\
(iii) There exist a constant $C_3 \in (0,\infty)$ and a function 
$G:\Re^k\rightarrow [0,\infty)$ with $E G(\hat{\udelta})<\infty$
 such that 
$$
|M_{2i}(\ux)| \leq m^{-1}G(\ux) \qmq{for all} \ux \in \Delta, 
$$
and 
$$
|M_{2i}^{(j)}(\udelta)| \leq C_3 m^{-1},$$
{for all} 
$j=1,\cdots,k; i=1,\cdots,m, m\geq 1.$
\end{enumerate}

We now briefly comment on the regularity condition S.
Condition S requires the functions $M_{1i}$ and $M_{2i}$
to be smooth, which typically holds under suitable smoothness conditions
on the parametric model (2.1) and (2.2). As mentioned
earlier, in most applications the function $M_{1i}$ is of the 
order $O(1)$ while $M_{2i}$ is of the 
order $O(m^{-1})$ as $\mti$.  Condition S requires that the 
partial derivatives of these functions also have the 
same orders, respectively. Condition S.3(iii) is a local 
Lipschitz condition of order $\eta\in (0,1]$ on 
 $M_{1i}$  and $M_{2i}$. This condition holds
with $\eta=1$ if  $M_{1i}$ is three-times continuously
differentiable  and $M_{2i}$ two-times continuously
differentiable  on a neighborhood of the true parameter value
$\udelta$.

Next, suppose that the bias and the variance matrix of 
the  given estimator $\hat{\udelta}$ are of the form:
\bea
\ub \equiv  E(\hat{\udelta}-\udelta)& = & \frac{\ua}{m}
+o(\frac 1m)\qmq{as}\mti\\
V\equiv Var (\hat{\udelta}) & = & \frac{\uSigma}{m}
+o(\frac 1m)\qmq{as}\mti.
\eea
Let $\hat{\ua}$ and $\hat{\uSigma}$ be estimators of 
the parameters $\ua$ and $\uSigma$ in 
(4.1) and (4.2) respectively,
such that for 
some $\eta \in (0,1]$, 
\bea
E\|\hat{\ua}-\ua\|^{1+\eta} & = & o(1)\qmq{as}\mti\\
E\|\hat{\uSigma}-\uSigma\|^{1+\eta} & = & o(1)\qmq{as}\mti.
\eea
Note that in the notation of Section 3, the quantities 
$\hat{\ub}$, $\hat{\ua}$, $\hat{\uV}$ and $\hat{\uSigma}$
are related as 
$\hat{\ub}=m^{-1}\hat{\ua}$ and 
$\hat{\uV}= m^{-1} \hat{\uSigma}$.

With this, we are now ready to state the main results of this section.
The first result shows that 
the   preliminary-titled-estimators $\bar{\udelta}_i$ converge
to the   true parameter  $\udelta$  
in probability  uniformly in $i=1,\ldots,m$,
and also that the first order
partial derivative  $M_{1i}^{(1)}(\hat{\udelta})$ falls
below the given threshold $(1+\log m)^{-2}$ with very small 
probability, uniformly in $i=1,\ldots,m$.\\[.2in]
{\bf  Theorem 1:}~ Suppose that (4.1)-(4.4) and condition S 
 hold. Then\\ 
(i) for any $\epsilon\in (0,\infty)$,
$$
\max_{1\leq i\leq m} 
	P\Big(\|\bar{\udelta}_i -\udelta\| >\epsilon\Big) = O(m^{-1})
\qmq{as}\mti.\\
$$
(ii) As $\mti$,  
$$\max_{1\leq i\leq m} P\Big( |M_{1i}^{(1)}(\hat{\udelta})| <
(1+\log m)^{-2}\Big) = O(m^{-1}).
$$

\noindent
{\bf Proof:} A  proof of the theorem is given in 
the Appendix.\\[.1in]

As a direct consequence of the above result, we get the following.
\\[.2in]
{\bf  Theorem 2:}~ Under the conditions of Theorem 1,
$$
\max_{1\leq i\leq m} 
	P\Big(~\tilde{\udelta}_i = \bar{\udelta}_i~ \Big) 
= 1- O(m^{-1})\qmq{as}\mti.\\
$$

\noindent
{\bf Proof:} A  proof of the theorem is given in 
the Appendix.\\[.1in]

Theorem 2 shows that uniformly in $i$, the titled estimator 
$\tilde{\udelta}_i$ coincides with 
the preliminary-titled-estimator $\bar{\udelta}_i$ with high 
probability when $m$ is  large. Thus, the typical
value of the titled estimator has a  correction term 
added to the first component of the given initial estimator 
$\hat{\udelta}$ (cf. (3.6)). The next  result shows  that this correction 
factor indeed reduces  the bias of the proposed MSPE 
estimator $\widehat{M_{i}(\udelta)}$ to order $o(m^{-1})$, 
as desired.\\[.2in]
{\bf  Theorem 3:}~ Suppose that (4.1)-(4.4) and condition S 
 hold. Further suppose that 
\be
\left\{\|\sqrt{m}(\hat{\udelta}-\udelta)\|^2\right\}_{m\geq 1}
\ee
is uniformly integrable. Then 
\be
\max_{1\leq i \leq m} 
|E\widehat{M_i(\udelta)}-M_i(\udelta)|=o(m^{-1})\qmq{as}\mti.
\ee

\noindent
{\bf Proof:} A proof of the theorem is given in the 
Appendix.\\[.1in]

\section{\sc Simulation Study}
\setcounter{section}{5}
\setcounter{equation}{0}

We conduct a small simulation study to check small sample 
performance of our  proposed MSPE estimator and compare 
it  with its competitors. In order to mimic a real 
life study, we consider the example in Battese, Harter and 
Fuller (1988) to  estimate the area under corn and soybeans
 for twelve counties of north-central Iowa. Originally, Battese 
{\it et al.} (1988) applied a nested error regression
 model. We consider here the area level version of their 
model for simplicity and we think that this is adequate 
for illustration purposes. Let $y_{ij}$ 
be the area under corn for $j$-th segment in $i$-th county 
and let  $\bar{X}_i$ be the (population) average number 
of pixels classified as corn  in the $i$-th county. We consider 
the area level model as 
\be
\bar{y}_i = \beta_0+\beta_1\bar{X}_i+u_i+e_i,\; i=1,\cdots, m
\ee
where $\bar{y}_i=\frac{1}{n_i}\sum_{j=1}^{n_i}y_{ij}=$ 
the sample average area under corn in the $i$-th county. 
Here, $u_i$'s are independently distributed with each
following  the  
$N(0,\sigma_u^2)$ distribution
and the $e_i$'s are independent with 
$e_i\sim N(0,D_i)$ for $i=1,\ldots,m$ where $D_i
=\frac{\sigma_e^2}{n_i}$. Further,  the  $u_i$'s and the $e_i$'s
are independent. In our simulation, we take $\beta_0=43.00, 
\beta_1=0.25, \sigma_u^2=140.00,  \sigma_e^2=147.00$. 
The $n_i$'s are as given in Battese {\it et al.} (1988) with 
$\min_{1\leq i\leq m}  {n_i}=1$,  ${\max}_{1\leq 
i\leq m}{n_i}=6$ and $m=12$. For the simulation study, 
we  generated
 $R=20,000$  sets of samples using model (5.1) and computed 
$\hat{\udelta}=(\hat{\beta}_0,\hat{\beta}_1,\hat{\sigma}_u)^T$ 
each time.

For estimating the bias and the variance of the estimator
vector $\hat{\udelta}$ used in the definition of the 
titled estimators $\tilde{\udelta}_i$'s, we employed a 
parametric bootstrap  method. For the sake of completeness, 
here  we briefly point out  the 
main steps of the bootstrap procedure.
%
%
\begin{itemize}
\item
Step (I):~  Generate independent random variables 
$\{e_i^*\}_{i=1}^m$ and $\{u_i^*\}_{i=1}^m$ with 
 $e_i^* \sim N(0,D_i)$ and 
$u_i^* \sim N(0,\hat{\sigma}_u^2)$.
\item
Step (II):~ 
 Define the bootstrap variables, $y_i^*
=\theta_i^*+e_i^*; ~ \theta_i^*=\ux_i^T\hat{\ubeta}
+u_i^*; i=1,\cdots, m.$
\item
Step (III):~ 
Define the  bootstrap version $\udelta^*$ of $\hat{\udelta}$ 
by replacing $y_1,\cdots, y_m$ in $\hat{\udelta}$ 
with $y_1*,\cdots, y_m^*.$ 
\end{itemize}
The bootstrap estimators of the bias and the variance 
matrix of $\hat{\udelta}$ are now given by 
\bea
\hat{\ub} & = & E_*\udelta^*-\hat{\udelta} \\ 
\hat{V} & = & E_*(\udelta^*-E_*\udelta^*)(\udelta^*-E_*\udelta^*)^T
\eea
where $E_*$ denote the conditional expectation given the data. 
In simulation, Steps (I)-(III) are repeated a large number of 
times and the average of the  bootstrap versions 
 $\udelta^*$'s gives the Monte-Carlo approximation to 
$E_*\udelta^*$ while the sample covariance matrix of 
the $\udelta^*$'s give the numerical value of the right side of
(5.3).

Next   for each of the three  MSPE estimators (namely,
the Prasad-Rao estimator
$\hat{M}_i^{PR}$, 
the Jiang, Lahiri and Wan estimator $\hat{M}_i^{JLW}$,
  and the proposed
estimator $\widehat{M_i(\udelta)}$)  of the  small area 
parameter $\theta_i$, we calculate 
the following measures:
\begin{itemize}
\item Relative bias with respect to the  empirical MSPE:
$$RB_i=\frac{E\{\widehat{MSPE(\hat{\theta}_i)}\}-SMSPE(\hat{\theta}_i)}
{SMSPE(\hat{\theta}_i)}, i=1\cdots,12$$
where, $ E\{\widehat{MSPE(\hat{\theta}_i)}=\frac 1R\sum_{r=1}^R 
\widehat{MSPE(\hat{\theta}_i)}^{(r)}$.

\item Empirical coefficient of variation:
$$CV_i=\frac{E^{\frac 12}\{\widehat{MSPE(\hat{\theta}_i)}
 -SMSPE(\hat{\theta}_i)\}^2}{SMSPE(\hat{\theta}_i)}, i=1,\cdots, 12 $$
where, $E\{\widehat{MSPE(\hat{\theta}_i)}
 -SMSPE(\hat{\theta}_i)\}^2=\frac 1R\sum_{r=1}^R\{\widehat
{MSPE(\hat{\theta}_i)}^{(r)}-SMPE(\hat{\theta}_i)\}^2$. 
\end{itemize}

Table 1 reports a summary result of the simulation study. The proposed 
estimator is denoted as `New' in the table. 
\begin{center}
Table 1: Summary of simulation study\\

	(a)~~Relative Bias

\begin{tabular}{lllllll}\\ 
    & min & $Q_1$ & median & mean & $Q_3$ & max \\
\\
PR  & -.164 & -.114 & -.055 & -.061 & -.025 & .048 \\
JLW & -.210 & -.142 & -.067 & -.095 & -.040 & -.022 \\
New & -.163 & -.113 & -.054 & -.060 & -.024 & .048
\end{tabular}

\vspace*{.2in}
(b)~~Empirical CV

\begin{tabular}{lllllll}\\ 
    & min & $Q_1$ & median & mean & $Q_3$ & max \\
\\
PR  & .010 & .033 & .055 & .074 & .114 & .164 \\
JLW & .082 & .094 & .120 & .132 & .151 & .212 \\
New & .009 & .034 & .054 & .074 & .113 & .163 
\end{tabular}

\end{center}

 From the above  table, it is clear that the proposed
estimator  and the  PR estimator performs at par 
and both perform  better than the jackknife-based 
 estimator, particularly in terms of the coefficient of 
variation. We should also mention that, in this 
simulation study, fortunately the jackknife 
method did not produce any negative MSPE estimates. 
This is perhaps due to the fact that 
the true parameter values are far 
away from the boundary of the parameter space. 

\section{\sc Conclusions}

In this paper, we described a new method of  bias 
correction for the naive `plug-in' estimator
of the MSPE of a function of the small area means
$h(\theta_i)$, $i=1,\ldots,m$. Unlike the existing methods
which may produce a negative estimate of the MSPE with positive
probability, the estimates of the MSPE produced by the  proposed
method   is always  nonnegative.  Theoretical properties 
of the method are investigated, which in particular 
show that the resulting estimator of the 
MSPE attains the same level of accuracy as the existing 
methods in correcting the bias of the naive  
 MSPE estimator. Further, the  numerical results presented in 
the paper shows that the proposed method performs at per 
with the  Prasad-Rao (1990) method, and has a slightly
better performance compared to the estimator based on
the jackknife method. A key difference of the new  method with 
the existing methods is that while the existing methods
apply various bias correction techniques to the MSPE 
function itself, the  new method reduces the  bias
 implicitly by  suitably tilting the value of argument
of the MSPE function.

\vspace*{1in}
\begin{center}
{\huge {\sc Appendix}}\\[.3in]
\end{center}
\renewcommand{\theequation}{A.\arabic{equation}}
\setcounter{equation}{0}
For a vector $\ux\in\Re^k$, let $\ux(j)$ denote the 
$j$th component of $\ux$, $j=1,\ldots,k$. Let $C, C(\cdot)$ 
 denote generic positive constants  that may depend on 
the argument(s) (if any) but not on $i=1,\ldots, m$  
or $m$. Also, unless explicitly specified, limits in order 
symbols are taken letting $\mti$.\\[.1in]

\noindent
{\bf Proof of Theorem 1:}~Since $M_{1i}(\udelta)\neq 0$,
by condition S.3.(ii), there exists $\epsilon_1, 
\epsilon_2\in (0,
\epsilon_0)$ such that  $|M_{1i}(\ux)|> \epsilon_2$
for all $\ux$ with $\|\ux -\udelta\| \leq \epsilon_1$. 
Hence, again by  condition S,
there exists a constant $C=C(\epsilon_2)\in (0,\infty)$
such that on the set $\{\|\hat{\udelta} -\udelta\|\leq 
\epsilon_1$,
$$
\|\hat{\udelta}_i -\hat{\udelta}\| \leq C\Big[\|\hat{\ub}\|
+\|\hat{V}\|\Big]
$$
uniformly in $i=1,\ldots,m$, $m\geq 1$. Hence, by Chebychev's
inequality, for any $\epsilon\in (0,\infty)$, 
\bea
&&\max_{i=1,\ldots,m} 
	P\Big(\|\bar{\udelta}_i -{\udelta}\| > \epsilon\Big)\nn\\
&\leq & P\Big(\|\hat{\udelta} -{\udelta}\| > \epsilon_1\Big)
	+ P\Big( C [\|\hat{\ub}\|
		+\|\hat{V}\|]> \epsilon\Big)\nn\\
&\leq & \epsilon_1^{-2} E\|\hat{\udelta} -{\udelta}\|^2
	+ \epsilon^{-1} CE\Big[\|\hat{\ub}\|
		+\|\hat{V}\|\Big]\nn\\
&=& O(m^{-1}).\nn
\eea
This proves part (i). For part (ii), note that 
\bea
&&\max_{i=1,\ldots,m} 
	P\Big(|M_{1i}^{(1)}(\hat{\udelta})| \leq (1+\log m)^{-2}\Big)\nn\\
&\leq & P\Big(|M_{1i}^{(1)}(\hat{\udelta}) -
	M_{1i}^{(1)}({\udelta})| > \epsilon_2/2\Big)\nn\\
&\leq & P\Big(\|\hat{\udelta} -{\udelta}\| > \epsilon_1\Big)
 + P\Big(C\|\hat{\udelta} -{\udelta}\|^{\gamma} > \epsilon_2/2\Big)
\nn\\
&=&  O(m^{-1}).\nn
\eea

\vspace*{.2in}

\noindent
{\bf Proof of Theorem 2:}~
Since $\udelta$ is an interior point of $\Delta$,
there exists a $\epsilon_3\in(0,\epsilon_0)$ 
such that $\{\ux : \|\udelta -\ux\|\leq \epsilon_3\}
\subset \Delta$.
Hence, by Theorem 1, 
\bea
&&\max_{i=1,\ldots,m} 
	P\Big(~\tilde{\udelta}_i \neq{\udelta}~\Big)\nn\\
&\leq& P(\hat{\udelta}\notin \Delta) +
\max_{i=1,\ldots,m} 
	P\Big(|M_{1i}^{(1)}(\hat{\udelta})| \leq (1+\log m)^{-2}\Big)\nn\\
&\leq & P\Big(\|\hat{\udelta} -{\udelta}\| > \epsilon_3\Big)
 + O(m^{-1})\nn\\
&=&  O(m^{-1}),\nn
\eea
where the last step follows by an application of Chebychev's inequality
as in the proof of Theorem 1 above. This proves Theorem 2.

\vspace*{.2in}

\noindent
{\bf Proof of Theorem 3:}~ By Taylor's expansion and condition S, on 
the set $\{\hat{\udelta} \in \N\}$, 
\be
M_{2i}(\hat{\udelta})=M_{2i}(\udelta)+
(\hat{\udelta}-\udelta)^T\triangledown M_{2i}(\udelta)+R_{1i}
\ee
where $\triangledown M_{2i}(.)$ is the $k \times 1 $ 
vector of first order
 partial derivatives of $M_{2i}$ and $R_{1i}$ is a 
remainder term. By condition S, $R_{1i}$ admits the bound 
\be
|R_{1i}| \leq \|\hat{\udelta}-\udelta\|\|\triangledown M_{2i}(\udelta)-
\triangledown M_{2i}(\udelta^0)\| \leq C m^{-1}\|\hat{\udelta}-
\udelta\|^{1+\gamma}
\ee
uniformly in $i=1,\cdots, m, m\geq 1$ where $\udelta^0$ is 
a point on the 
line joining $\hat{\udelta}$ and $\udelta$, so that 
$\|\udelta^0-\udelta\| \leq \|\hat{\udelta}-\udelta\|$. 
Hence, by (3.1), (3.2), (A.1), (A.2) and the dominated 
convergence theorem (DCT), 

\bea
&&\max_{i\leq i \leq m} |E M_{2i}(\hat{\udelta})-M_{2i}(\udelta)|\nn\\
&\leq & \max_{1\leq i \leq m} 
|E \{M_{2i}(\hat{\udelta})-
	M_{2i}(\udelta)\}\amalg(\hat{\udelta}\in \N)| 
+ \max_{1\leq i \leq m}
 E\{ M_{2i}(\hat{\udelta})+M_{2i}(\udelta)\}\amalg (\hat{\udelta} 
\notin \N)
 \nn \\
&\leq & \max_{1 \leq i \leq m} \left[ \|E(\hat{\udelta}-
\udelta)\amalg(\hat{\udelta}\in \N)\| \cdot
\|\triangledown M_{2i}(\udelta)\| 
+E|R_{1i}|\amalg(\hat{\udelta}\in \N)\right] \nn \\
 &&  + m^{-1}E\{G(\hat{\udelta})+G(\udelta)\}\amalg(\hat{\udelta}
\notin \N)\nn\\
 &\leq &  \max_{1 \leq i \leq m} \left[\left\{ 
\|E(\hat{\udelta}-\udelta)\|
+E\|\hat{\udelta}-\udelta\|\amalg(\hat{\udelta}\notin \N)\right\}
\|\triangledown M_{2i}(\udelta)\| 
+C m^{-1}E\|\hat{\udelta}-\udelta
\|^{1+\gamma}\right] \nn \\
 && +Cm^{-1}\left[P(\hat{\udelta}\notin \N)
+EG(\hat{\udelta})\amalg(\hat{\udelta}\notin \N)\right] 
\nn\\
&\leq & C\left[ m^{-2}+m^{-1}(E\|\hat{\udelta}-\udelta\|^2)^{\frac 12} 
\left(P(\hat{\udelta} \notin \N)\right)^{\frac 12} +m^{-1}\left(E\|
\hat{\udelta}-\udelta\|^2\right)^{\frac{1+\gamma}{2}}\right] +o(m^{-1})\nn \\
&= & o(m^{-1}),
\eea
as $P(\hat{\udelta}\notin \N)
 \leq \epsilon_0^{-2}E\|\hat{\udelta}-\udelta\|^2= O (m^{-1})$. 
Without loss of generality, suppose that $\epsilon_0$ (in the definition of 
 $\N$) is small enough so that for some $C\in (0,\infty),~
 \sup\{|M_{1i}^{(1)}(\ux)|^{-1}: \ux \in \N, j,l=1,\cdots, k;
 i=1,\cdots,m\} <C$.  Let 
\bea
\hat{L}_i &=& \sum_{j=1}^k M_{1i}^{(j)}(\hat{\udelta})\hat{\ub}(j)+
	\sum_{j=1}^k\sum_{l=1}^kc(j,l)M_{1i}^{(j,l)}(\hat{\udelta})
			\hat{V}(j,l)\nn \\
\tilde{L}_i &=& \sum_{j=1}^k M_{1i}^{(j)}({\udelta})\hat{\ub}(j)+
	\sum_{j=1}^k\sum_{l=1}^kc(j,l)M_{1i}^{(j,l)}({\udelta})
			\hat{V}(j,l),\nn 
\eea
$i=1,\ldots,m$, where $c(j,l) =1/2$ for $j\neq l$ and
 $c(j,l) =1$ for $j= l$.
Then by Taylor's expansion, it follows that 
there exists a constants $C \in (0, \infty)$ (not depending on $i$) 
such that on the set $\{\hat{\udelta}\in \N\}$, 
\bea
\hat{L}_i & = &  \tilde{L}_i +R_{2i}, \mbox{ say}
\eea
and 
\be
\bar{\udelta}_i 
= \hat{\udelta}-\frac{\tilde{L}_i}{M_{1i}^{(1)}(\udelta)}\ue_1+R_{3i}\ue_1
\ee
for $i=1,\cdots,m$ where 
$\max_{1 \leq i \leq m} |R_{2i}| \leq C\Big\{\|\hat{\ub}\|.
\|\hat{\udelta}-\udelta\|+\|\hat{\udelta}-\udelta\|^{\gamma}\|\hat{\uV}\|
\Big\}$ and for all $i=1,\cdots, m,$
\be
|R_{3i}|\leq C\Big\{|\hat{L}_i|.\|\hat{\udelta}-\udelta\|+|R_{2i}|\Big\}.
\ee

By similar arguments, on the set $A_{1i}
\equiv \{\hat{\udelta}\in N\}
\cap\{\bar{\udelta}_i\in \Delta\}$, we may write
\bea
&&\sum_{j=1}^k\sum_{l=1}^kc(j,l)M_{1i}^{(jl)}(u\bar{\udelta}_i+
(1-u)\hat{\udelta})(\bar{\udelta}_i(j)-\udelta(j))
(\bar{\udelta}_i(l)-\udelta(l))\nn\\
&=&\sum_{j=1}^k\sum_{l=1}^kc(j,l)M_{1i}^{(jl)}
(\udelta)(\hat{\udelta}_i(j)-\udelta(j))
(\hat{\udelta}_i(l)-\udelta(l))+R_{4i}(u)\nn
\eea
$i=1,\cdots,m; u\in [0,1]$, where 
\be
\sup_{u\in [0,1]} \max_{1\leq i \leq m} |R_{4i}(u)| \leq 
C\left[\|\bar{\udelta}_i-\udelta\|^{2+\gamma}
+\|\bar{\udelta}_i-\hat{\udelta}\|\cdot \|\hat{\udelta}-\udelta\| 
+\|\bar{\udelta}_i-\hat{\udelta}\|^2\right]
\ee
for some $C\in (0,\infty).$

On the set $A_{2i}=A_{1i}\cap \{\bar{\udelta}_i\in \N\}
=\{\hat{\udelta}\in \N\}
\cap\{\bar{\udelta}_i\in \N\}$, by Taylor's expansion, there exists a point 
$\udelta_i^*$ on the line joining $\bar{\udelta}_i$ and $\udelta$ such that 
\bea
&&M_{1i}(\bar{\udelta}_i)-M_{1i}(\udelta)\nn\\
&=& \sum_{j=1}^kM_{1i}^{(j)}(\udelta)\{\bar{\udelta}_i(j)-\udelta(j)\}
+\sum_{j=1}^k\sum_{l=1}^k c(j,l)M_{1i}^{(jl)}(\udelta_i^*)
\{\bar{\udelta}_i(j)-\udelta(j)\}\{\bar{\udelta}_i(l)-\udelta(l)\} \nn \\
& =&  \sum_{j=1}^kM_{1i}^{(j)}(\udelta)\left(\hat{\udelta}(j)-\udelta(j)\right)
+M_{1i}^{(1)}(\udelta)\left(-\frac{\tilde{L}_i}{M_{1i}^{(1)}(\udelta)}+
R_{3i}\right) \nn \\
 && +\sum_{j=1}^k\sum_{l=1}^k c(j,l)M_{1i}^{(jl)}(\udelta)\{\hat{\udelta}(j)-
\udelta(j)\}\{\hat{\udelta}(l)-\udelta(l)\}+R_{4i}^* \nn\\ 
&=&  \sum_{j=1}^kM_{1i}^{(j)}(\udelta)
	\left\{\Big(\hat{\udelta}(j)-\udelta(j)\Big) 
		-\hat{\ub}(j)\right\}\nn\\
&&\hspace{.3in}+\sum_{j=1}^k\sum_{l=1}^k c(j,l)
	M_{1i}^{(j,l)}(\udelta)\left\{
		\Big(\hat{\udelta}(j)-\udelta(j)\Big)
			\Big(\hat{\udelta}(l)-\udelta(l)\Big)
				-\hat{V}(j,l)\right\}\nn\\
&&\hspace{1in}+M_{1i}^{(1)}(\udelta)R_{3i}+R_{4i}^* \nn \\ 
& \equiv&  Q_{1i}+M_{1i}^{(1)}(\udelta)R_{3i}+R_{4i}^*, \mbox{ say }
\eea
where $R_{4i}^*=R_{4i}(u)$ with the $u$ corresponding to $\udelta_i^*$.

Hence for $i=1,\cdots, m$, 
 with $A_{3i}
=\{|M_{1i}^{(1)}(\hat{\udelta})|^{-1}\leq (1+\log m)^2\},$ 
\bea
&&M_{1i}(\tilde{\udelta})-M_{1i}(\udelta)\nn\\
& =&  [M_{1i}(\bar{\udelta}_i)-M_{1i}(\udelta)]
\amalg\Big(\left\{\bar{\udelta}_i\in 
		\Delta\right\}\cap A_{3i}\Big)
+[M_{1i}(\hat{\udelta})-M_{1i}(\udelta)]
\amalg\Big(\left\{\bar{\udelta}_i\notin 
		\Delta\right\}\cup A_{3i}^c\Big)\nn \\ 
&=&  [M_{1i}(\bar{\udelta}_i)-M_{1i}(\udelta)]\Big\{\amalg(A_{2i}) 
+\amalg (\bar{\delta}_i\in \Delta)-
	\amalg(A_{2i})\Big\}\amalg(A_{3i}) \nn\\
&&\hspace{.3in}+
[M_{1i}(\hat{\udelta})-M_{1i}(\udelta)]
	\amalg\Big(\left\{\bar{\udelta}_i\notin 
		\Delta\right\}\cup A_{3i}^c\Big) \nn \\ 
& \equiv&  \Big[Q_{1i}+M_{1i}^{(1)}(\udelta)R_{3i}+R_{4i}^*\Big]
\amalg(A_{2i}\cap A_{3i}) +R_{5i}, \mbox{ say } \nn \\ 
& \equiv &  Q_{1i}+R_{6i}, \mbox{say,} 
\eea
where $|R_{6i}| \leq |R_{5i}|+|R_{3i}+R_{4i}^*|\amalg(A_{2i})
+|Q_{1i}|\amalg(A_{2i}^c\cap A_{3i}^c)$ and 
\bea
|R_{5i}| & \leq & |M_{1i}(\bar{\udelta}_i)-M(\udelta)|
	\cdot |\amalg(\bar{\udelta}_i\in \Delta)
		-\amalg(A_{2i})|\amalg(A_{3i}) \nn \\
 & & \hspace{.3in}+ |M_{1i}(\hat{\udelta}-M_{1i}(\udelta)|
\amalg(\{\bar{\udelta}_i \notin \N\}\cup A_{3i}^c) \nn \\ 
 & \equiv & R_{51i}, \mbox{ say. } \nn
\eea
Note that by definition, 
\bea 
&&\mid \amalg (\bar{\udelta}_i\in \Delta) -\amalg (A_{2i})\mid
\nn\\
& \leq &  \amalg(\bar{\udelta}_i \in \Delta) \amalg(A_{2i}^c) 
+\amalg(\bar{\udelta}_i \notin \Delta) \amalg (A_{2i}) \nn\\ 
& \leq & \{ \amalg (\hat{\udelta} \notin \N) + \amalg(\bar{\udelta}_i 
\in \Delta \setminus \N)\} + \amalg (\emptyset). \nn \\
\eea
Hence, with $A_{4i}^c=\{\bar{\udelta}_i\notin \N\}\cap A_{3i}$, 
\bea
R_{51i} & \leq & \mid M_{1i}(\bar{\udelta}_i)-M_{1i}(\hat{\udelta})\mid
 \amalg (A_{3i})\{\amalg(\hat{\udelta} \notin \N) +\amalg(\bar{\udelta}_i 
\in \Delta \setminus\N)\} \nn \\ 
 & & + 2 \mid M_{1i}(\hat{\udelta})
	-M_{1i}(\udelta)\mid
		\Big\{\amalg(\hat{\udelta} \notin N)
		+\amalg\Big(\{\bar{\udelta}_i\notin N\}\cap A_{3i}\Big)
			+\amalg(A_{3i}^c)\Big\}  \nn \\
& \leq & C \|\bar{\udelta}_i-\hat{\udelta}\|\amalg(A_{3i})
	\Big\{\amalg(\hat{\udelta}\notin N) + \amalg(\bar{\udelta}_i\in 
		\Delta \setminus \N)\Big\} \nn\\
& & + C\|\hat{\udelta}-\udelta\|\Big\{\amalg(\hat{\udelta}\notin \N) 
	+ \amalg (A_{4i}^c) +\amalg(A_{3i}^c)\Big\}  \nn \\
& \leq & C \cdot (\log m)^2\{\|\hat{\ub}+\|\hat{\uV}\|\}
	\{\amalg(\hat{\udelta} \notin \N)+\amalg(A_{4i}^c)\} \nn\\
 & & + C \cdot \|\hat{\udelta}-\udelta\|
		\Big\{\amalg(\hat{\udelta} \notin \N) 
	+\amalg(A_{4i}^c)+\amalg(A_{3i}^c)\Big\}. 
\eea
By condition S, there exist $C \in (0, \infty)$ 
and $\epsilon_1 \in (0, \frac{\epsilon_0}{2})$ such that 
\bea
A_{4i}^c & \subset & 
	\{\|\hat{\udelta}-\udelta\| > \frac{\epsilon_0}{2}\}\cup
	\{\|\bar{\udelta}_i-\hat{\udelta}\| > \frac{\epsilon_0}{2}\}
		\nn \\ 
 & \subset & \{\|\hat{\udelta}-\udelta\|>\epsilon_1\} 
	\cup \{(\log m )^2 (\|\hat{\ub}\|+\|\hat{\uV}\|)>C\} 
\eea
and $A_{3i}^c \subset \{\|\hat{\udelta}-\udelta\| >\epsilon_1\}$
 for all $i=1,\cdots,m, m\geq 1$. Hence, it follows that 
\bea
R_{51i} & \leq & C \cdot (\log m)^2 \{ \| {\hat{\ub}}\| +\|{\hat{\uV}}\|\} 
	\left[	\amalg(\|\hat{\udelta}-\udelta\|>\epsilon_1) 
		+\amalg\Big([\log m]^2(\|\hat{\ub}\|
			+\|\hat{\uV}\|)>C\Big)\right] \nn \\ 
 && + C\cdot  \|\hat{\udelta}-\udelta\|\left[
\amalg (\|\hat{\udelta}-\udelta\|>\epsilon_1) 
+\amalg\Big([\log m]^2(\|\hat{\ub}\|+\|\hat{\uV}\|)>C\Big)\right]  
\eea
for all $i=1,\cdots,m, m \geq 1$. Let $W_1=(\|\hat{\ua}\|
+\|\hat{\uSigma}\|)$. Note that by uniform integrability of
 $\{({\sqrt{m}}\|\hat{\udelta}-\udelta\|)^2\}_{m \geq 1}$ 
and the fact that $E\mid W_1\mid^{1+\eta} = O (1)$, 
\bea
&&\max_{1 \leq i \leq m} E (R_{51i})\nn\\
 & \leq & C m^{-1} (\log m)^2
\left[ \Big(E\mid W_1\mid^{1+\eta}\Big)^{\frac{1}{1+\eta}}
\Big(P (\|\hat{\udelta}-\udelta\|>\epsilon_1\Big)
^{\frac{\eta}{1+\eta}}
+E \mid| W_1\mid^{1+\eta}\{m^{-1}(\log m)^2\}^{\eta}\right]\nn  \\ 
 && + C\left[\epsilon_1^{-1}E\|\hat{\udelta}-\udelta\|^2
\amalg(\|\hat{\udelta}-\udelta\|>\epsilon_1)
+\left(E\|\hat{\udelta}-
\udelta\|^2\right)^{1/2}\Big\{P(m^{-1}(\log m)^2 
| W_1 | >C)\Big\}^{\frac 12} \right] \nn \\ 
& = & o(m^{-1}) \qmq{as}\mti.
\eea 
This completes the proof of Theorem 3.\\[.2in]

\noindent
{\bf Acknowledgment:} The research is partially supported by NSF 
grant number  DMS 0306574 and SES 0318184.\\[.3in]

\begin{center}

{\sc References}
\end{center}

\begin{description}
\item{\sc Battese, G.E., Harter, R.M., \& Fuller, W.A.} (1988). 
An error component model for prediction of county crop areas using 
survey and satellite data, {\em J. Am. Statist. Assoc.} {\bf 83}, 28-36. 

\item{\sc Bell,W.} (2002). Discussion with `` Jackknife in the 
Fay-Herriot model with an application'', {\em Proceeding of the seminar on 
funding opportunity in survey research}, 98-104.

\item{\sc Datta, G.S. \& Lahiri, P.} (2000). A unified measure of 
uncertainty of estimated best linear unbiased predictors in small area 
estimation problems, {\em Statistica Sinica} {\bf 10}, 623-27.

\item{\sc Efron, B. } (1979).
     Bootstrap methods: Another look at the jackknife, 
{\em  Ann.  Statist.} {\bf  7}, 1-26.  

\item{\sc Jiang, J., Lahiri, P., \& Wan, S-M.} (2002). A unified 
jackknife theory for empirical best prediction with M-estimation, 
{\em Ann. Statist.} {\bf 30}, 1782-810.

\item{ \sc Lahiri, P. \& Rao, J.N.K.} (1995). Robust estimation of mean 
squared error of small area estimators, {\em J. Am. Statist. Assoc.} {\bf 82},
 758-66. 
 
\item{\sc Prasad, N.G.N. \& Rao, J.N.K.} (1990). The estimation of the 
mean squared error of small area estimators. {\em J.Am. Statist. Assoc.} 
{\bf 68}, 67-72. 

\item{\sc Rao, J.N.K.} (2003). {\em Small Area Estimation}, Wiley.
\end{description}

\end{document}